\documentclass[11pt]{article}

\usepackage{amsmath,amsthm,amssymb}

\renewcommand{\k}{\ensuremath{\Bbbk}} 
\newcommand{\double}[3]{{#1}{\backslash}{ #2}{/}{#3}} 

\newcommand{\sdouble}[3]{\size{\double{#1}{#2}{#3}}} 
\newcommand{\rightco}[2]{#1\backslash #2} 
\newcommand{\st}{{\ \mid \ }} 
\newcommand{\C}[1]{C_{#1}^{}} 
\newcommand{\size}[1]{\left|#1\right|}
\newcommand{\rtsys}{\Phi} 
\newcommand{\rtbase}{\Delta} 

\newcommand{\inv}{{}^{-1}}
\newcommand{\sizefrac}[2]{\frac{\size{#1}}{\size{#2}}}
\newcommand{\dsize}[1]{\left|#1\right|}
\renewcommand{\tilde}{\widetilde}
\renewcommand{\hat}{\widehat}
\renewcommand{\phi}{\varphi}
\renewcommand{\th}{\ensuremath{{}^{\text{th}}}}

\newenvironment{map_machine}[5]
   {
    \renewcommand{\arraycolsep}{1pt}
    \begin{array}[#1]{rcl}
   #2 & {-}\!\!{-}\!\!{\longrightarrow} & #3\\ 
   #4 & \mapstochar\kern-.25em{-}\!\!{\longrightarrow} & #5
   }{\end{array}}
\newcommand{\map}[5][t]{\begin{map_machine}{#1}{#2}{#3}%
   {#4}{#5}\end{map_machine}}

\newtheorem{theorem}{Theorem}[section]
\newtheorem{lemma}[theorem]{Lemma}
\newtheorem{cor}[theorem]{Corollary}
\newtheorem{prop}[theorem]{Proposition}
\theoremstyle{definition}
\newtheorem{remark}[theorem]{Remark}

\newtheorem{innerproof}[theorem]{Proof}
\newenvironment{proofof}
        {\begin{innerproof}}{\qed\end{innerproof}}

\newcommand{\define}[1]{\textbf{#1}}

\newcommand{\compress}
{\setlength{\partopsep}{0in}
      \setlength{\itemsep}{0in}
      \setlength{\topsep}{0in}
      \setlength{\parsep}{0in}
      \setlength{\parskip}{0in}}
\newenvironment{theorem_count}
     {\compress
      
      \begin{enumerate}
        \compress}
   {\end{enumerate}}

\DeclareMathOperator{\Aut}{Aut}
\DeclareMathOperator{\charistic}{char}

\DeclareMathOperator{\SL}{SL}
\DeclareMathOperator{\GL}{GL}
\newcommand{\F}{\ensuremath{\mathbb{F}}}
\newcommand{\Z}{\mathbb{Z}}
\newcommand{\E}{\mathbb{E}}
\DeclareMathOperator{\N}{\mathbb{N}}

\newbox{\myaddress}
\savebox{\myaddress}{{\rm \normalsize 334 Hill Center, Rutgers 
University, Piscataway NJ, 08854}}
\newbox{\myemail}
\savebox{\myemail}{\tt \normalsize duck@math.rutgers.edu}

\author{W. Ethan Duckworth\\
\usebox{\myaddress}\\
\usebox{\myemail}}

\title{Infiniteness of Double Coset 
Collections in Algebraic Groups}

\date{January 10, 2003.}

\begin{document}

\maketitle

\begin{abstract}
Let $G$ be a linear algebraic group defined over an algebraically
closed field.  The double coset question addressed in this paper is
the following: Given closed subgroups $X$ and $P$, is $\double XGP$
finite or infinite?  We limit ourselves to the case where $X$ is
maximal rank and reductive and $P$ parabolic.  This paper presents a
criterion for infiniteness which involves only dimensions of
centralizers of semisimple elements.  This result is then applied to
finish the classification of those $X$ which are spherical.  Finally,
excluding a case in $F_4$, we show that if $\double XGP$ is finite
then $X$ is spherical or the Levi factor of $P$ is spherical.  This
implies that it is rare for $\double XGP$ to be finite.  The primary
method is to descend to calculations at the finite group level and 
then to use elementary character theory.
\bigskip

\end{abstract}

\section{Introduction}
All \label{statement_of_results} algebraic groups in this paper are
linear algebraic groups defined over an algebraically closed field and
all subgroups are assumed to be closed.  Given an algebraic group $G$
we wish to classify those subgroups $X$ and $P$ such that $\double
XGP$ is finite.  One has that $\double XGP$ is finite if and only if
the $G$-orbit $G/P$ splits into finitely many $X$-orbits.  This
viewpoint makes a complete classification of all finite double coset
collections appear unlikely in the near future.  For this reason, we
will restrict which subgroups we allow.  We will generally assume $G$
is a reductive (or simple) algebraic group, that $P$ is a parabolic
subgroup and that $X$ is maximal rank and reductive.  This paper is
concerned with proving that a double coset collection is infinite.  We
intend to establish finiteness results in a later paper.

We will state the main results of the paper first with brief 
indications of how these results relate to earlier work in the 
field.  This is followed by a lengthier description of some of this 
earlier work.  


The first result provides a powerful criterion for establishing that
$\double XGP$ is infinite.  If $H$ is a group and $g\in H$ we write
$H_g$ for the centralizer of $g$ in $H$.

\begin{theorem}[Dimension Criterion]
Let \label{dimension_criterion} $G$ be a reductive algebraic group,
$X$ and $P$ subgroups of $G$ with $X$ maximal rank and $P$ parabolic. 
Let $L$ be a Levi factor of $P$ and let $s\in X\cap L$ be a semisimple
element.  If $\dim Z(G_s) + \dim G_s > \dim X_s + \dim P_s$
(equivalently, if $\dim Z(G_s) + \frac 12 \dim G_s - \dim X_s - \frac
12 \dim L_s >0$), then $\double {X_s}{G_s}{P_s}$ and $\double XGP$ are
infinite.
\end{theorem}

\paragraph{Classification of maximal rank reductive spherical subgroups.}
For \label{para_classification_of_spherical} the next result, suppose
$G$ is a simple algebraic group.  The first application of the
dimension criterion is to finish the classification of maximal rank
reductive spherical subgroups.  The subgroup $X$ is \define{spherical}
if a Borel subgroup $B$ has a dense orbit upon $\rightco XG$.  Brion
\cite{brion} and Vinberg \cite{vinberg} independently showed that $X$
is spherical if and only if $B$ has a finite number of orbits upon
$\rightco XG$, or, equivalently, $\double XGB$ is finite.  A maximal
rank reductive subgroup is \define{generic} if a subgroup of the same
type exists in all characteristics.  The work of Kr\"amer
\cite{kramer}, Brundan \cite{brundan} and Lawther \cite{lawther} has
produced a list of subgroups which are spherical in all
characteristics.  The generic maximal rank reductive subgroups on this
list are given in table \ref{spherical_subgroups} where we use the
following conventions.  We treat $A_0$ and $B_0$ as trivial groups and
$D_1$ as a $1$-dimensional torus.  We list only the Lie type of each
group, as the property of being spherical is not affected by which
representative of an isogeny class is used.  The notation $T_i$ refers
to an $i$ dimensional torus, central in $X$.  Finally, in the group
$A_1\tilde A_1\le G_2$ the tilde $\ \tilde{}\ $ is used to signify a
subgroup with short roots (we don't use this notation for the other
groups as there is no ambiguity).  Finally, in attempting to classify
the maximal rank reductive spherical subgroups, it suffices to
classify only the generic cases as the others arise from isogenies or
graph automorphisms, which preserve the property of being spherical
(see Lemma \ref{lemma_invariants_double_coset}).

%

\begin{table*}
\caption{Generic Maximal Rank Reductive Spherical Subgroups}
\renewcommand{\arraycolsep}{.2em}
$$\begin{array}{rclrcl}
\hline
X &\le&   G  & \qquad X & \le&  G\qquad \\[.5\jot]\hline\hline
\rule{3ex}{0in} 
A_nA_mT_1  & \le &  A_{n+m+1} & E_6T_1 & \le &   E_7 \\
B_nD_m     & \le &  B_{n+m}   &      A_7    & \le &   E_7  \\ 
A_{n-1}T_1 & \le &  B_n       &  A_1D_6 & \le &   E_7  \\
C_nC_m     & \le &  C_{n+m}   &  A_1E_7 & \le &   E_8\\    
C_{n-1}T_1 & \le &  C_n       &  D_8    & \le &   E_8\\     
A_{n-1}T_1 & \le &  C_n       &  A_1C_3 & \le &   F_4\\      
D_nD_m     & \le &  D_{n+m}   &  B_4    & \le &   F_4\\       
A_{n-1}T_1 & \le &  D_n       &  A_2   & \le  &  G_2\\       
D_5T_1     & \le &   E_6      &  A_1 \tilde A_1 & \le &   G_2\\
A_1A_5     & \le &   E_6\\                            
\hline
\end{array}$$
\label{spherical_subgroups}
\end{table*}

To prove that table \ref{spherical_subgroups} is complete, we
introduce the following root-theoretic property which was inspired by
Lawther's anti-open property (see Theorem \ref{lawthers_theorem}
below).  Let $X$ be a maximal rank reductive subgroup of $G$.  We
abbreviate the phrase ``maximal rank reductive'' with MRR. Fix a
maximal torus of $X$ and define $\rtsys(X)$ and $\rtsys(G)$ with
respect to this maximal torus.

We say $X$ has an \define{$R$-complement} if there exists a closed root
subsystem $R\le \rtsys(G)-\rtsys(X)$.  This is equivalent to the
existence of a generic MRR subgroup $K\le G$ with $\rtsys(K)= R$ and
$K\cap X$ a maximal torus.  The adjective ``long'' or ``short'' may be
applied if $R$ has only long or short roots.

\begin{theorem}
Let \label{completeness_of_table} $G$ be a simple algebraic group and
$X$ a generic MRR subgroup.  The following are equivalent:
\begin{theorem_count}
\item $X$ is spherical. 

\item $X$ appears in table \ref{spherical_subgroups}.

%
\item $X$ has no $A_2$ or $B_2$ complement.
\end{theorem_count}
\end{theorem}

In this paper we show that (i) $\implies$ (iii) $\implies$ (ii) (more
precisely, we show $\neg$ (ii) $\implies \neg$ (iii) $\implies \neg$
(i)).  The other implications are due to Brundan \cite{brundan} and
Lawther (see  Theorem \ref{lawthers_theorem} below).

Theorem \ref{completeness_of_table} applies to groups acting on the
full flag variety $G/B$, where $B$ is a Borel subgroup.  Using the
dimension criterion, we now obtain more general infiniteness results
where $P$ is a parabolic.  Since table \ref{spherical_subgroups}
contains relatively few subgroups, the following theorem places great
restrictions upon $X$ and $P$ for $\double XGP$ to be finite.  An
\define{end node parabolic} is conjugate to a standard parabolic
obtained by crossing off exactly one of the end nodes in the Dynkin
diagram of $G$.

\begin{theorem}[Spherical $X$ or Spherical $L$]
Let \label{sph_or_sph_lev} $G$ be a simple algebraic group, $X$ a MRR
subgroup, $P$ a parabolic subgroup with Levi factor $L$.  If $G=F_4$
suppose $P$ is not an end node parabolic.  If $\double XGP$ is finite
then $X$ is spherical or $L$ is spherical.
\end{theorem}

The extra restrictions placed upon $P$ when $G=F_4$ are necessary.  In
a later paper we will show that $\double {L_1}{F_4}{P_4}$ and $\double
{L_4}{F_4}{P_1}$ are finite (where $P_i$ is conjugate to the standard
parabolic obtained by crossing off the $i$\th\ node of the Dynkin
diagram of $F_4$, and $L_i$ is its Levi factor).

\begin{cor} 
If \label{spherical_or_maximal} $\double XGP$ is finite and $P$ is not
maximal then $X$ is spherical.
\end{cor}

\begin{remark} The theorem and the corollary give a surprisingly
strong dichotomy for MRR subgroups with respect to the double coset
problem.  Either they are spherical, or they will have an
infinite number of orbits on almost all flag varieties.  For instance,
$A_1A_5$ is spherical in $E_6$, but $X=T_1A_5$ will have an infinite
number of orbits on all flag varieties $E_6/P$ except, possibly, if
$P$ is an end node parabolic.  As another example, suppose one could
show that a MRR subgroup $X$ in $\GL(V)$ has a finite number of orbits
on flags consisting of one and two dimensional subspaces.  Then $X$
has a finite number of orbits on full flags, i.e. upon $G/B$ where $B$
is a Borel subgroup.  
\end{remark}


\section{Some History of the Problem}

\paragraph{Spherical subgroups.}
The \label{spherical_subgroup_background} following result is a
starting point for classifying reductive spherical subgroups.  Let $G$
be an algebraic group and let $\C G(\tau)$ be the centralizer of the
involution $\tau\in \Aut (G)$.

\begin{theorem}[Matsuki \cite{matsuki}, Springer \cite{springer}]
Let \label{centralizer_involution} $G$ be a reductive group.  If the 
characteristic of the underlying field is not $2$, then the 
centralizer of an involution in $G$ is spherical.  
\end{theorem}


This result was proven by Matsuki in characteristic $0$ and extended
to all characteristics except $p=2$ by Springer.  Seitz \cite{seitz}
has given an alternative proof and Lawther \cite{lawther} has given a
proof for positive characteristics which also shows that the same
subgroups are spherical in characteristic 2.  (Note, by ``same''
subgroups, we mean, by abuse of language, subgroups of the same type.) 
We will have occasion later to refer to Lawther's version of this
theorem.  In the following theorem let $\alpha_1$, \dots, $\alpha_n$
be the simple roots of the simple algebraic group $G$.  Write the high
root of $G$ as $\sum \lambda_i \alpha_i$.  Recall that a generic
maximal rank reductive subgroup of $G$ can be produced via the
Borel-de Siebenthal algorithm \cite{borel_desiebenthal}.  This
algorithm takes the Dynkin diagram, or the extended Dynkin diagram, of
$G$, crosses off some number of nodes, and recursively applies the
same procedure to the new diagrams which have been produced.

\begin{theorem}[Lawther \cite{lawther}]
Let 
\label{lawthers_theorem}
$G$ be a simple algebraic group with high root $\sum \lambda_i
\alpha_i$ as just described.  Let $X$ be a generic maximal rank
reductive subgroup.  The following conditions are equivalent and
define \emph{anti-open}.  If $X$ is anti-open then it is spherical.
\begin{theorem_count}
\item There do not exist $\alpha, \beta, \alpha+\beta \in
\rtsys(G)-\rtsys(X)$.  

\item $X$ can be produced from the Borel-de Siebenthal algorithm by
removing a single node $\alpha_i$ with $\lambda_i=1$ or by extending the
diagram for $G$ once and removing a single node $\alpha_i$ with
$\lambda_i=2$.

\item If the characteristic of the underlying field is not $2$, then $X$
is the centralizer of an involution.

\item $(G,X)$ appears in table \ref{spherical_subgroups} 
but $(G,X) \not \in \{ (B_n, A_{n-1}T_1)$, $(C_n, C_{n-1}T_1)$, $(G_2,A_2)\}$.

\end{theorem_count}
\end{theorem}

\begin{remark}
If $X$ is anti-open then $X$ has no $A_2$ or $B_2$ complement. 
Theorem \ref{lawthers_theorem} and Theorem \ref{completeness_of_table}
show that the converse is often, but not always, true.  To be precise,
the converse is true unless $(G,\ X) \in \{(B_n,\ A_{n-1}T_1)$, $(C_n,
C_{n-1}T_1)$, $(G_2,\ A_2)\}$.  For example, let $G=B_2$, fix a
maximal torus $T$ and label the Dynkin diagram of $G$ with $\alpha_1$
and $\alpha_2$ where $\alpha_1$ is long.  The group $X$, generated by
$T$ and the root groups corresponding to the high root and its
negative, is of type $A_1T_1$.  Then $\alpha_1$, $\alpha_2$, $\alpha_1
+ \alpha_2 \in \rtsys(G)-\rtsys(X)$ but $\alpha_1+2\alpha_2\in
\rtsys(X)$.  So $X$ does not have an $A_2$ or $B_2$ complement, but it
is not anti-open.
\end{remark}

If the characteristic $p$ is not $2$, then $\C G(\tau)$, the
centralizer of an involution, forms a reductive group.  The reductive
spherical subgroups of the simple algebraic groups have been
classified in characteristic 0 by Kr\"amer \cite{kramer} and in fact
most of them are centralizers of involutions.  Brundan \cite{brundan}
gave a reduction mod $p$ argument to show that certain groups in
Kr\"amer's list are spherical in all characteristics, but this
argument could not be applied to nine of the cases.  These nine
exceptions were all groups which would be the centralizer of an
involution when $p$ was not $2$.  Lawther \cite{lawther} introduced
the concept of anti-open to prove that these remaining cases were
spherical when $p=2$.  He then proceeded to work at the level of
finite groups using character theory.  See sections
\ref{section_reduction_to_finite_groups} and
\ref{section_dimension_criterion} below for some similar arguments.

In the work just described, the attention is focused on the group $X$. 
Related work has been done where the attention is focused on the
action of a Borel subgroup $B$ on varieties or modules.  For example,
the homogeneous space $\rightco {\C G(\tau)}G$ is a \emph{symmetric}
variety, and these have been studied by a wide variety of authors in
various contexts of representation theory.  See, for example, Helminck
\cite{helminck2} for how to compute the orbits of $B$ upon $\rightco
{\C G(\tau)}G$ and for a review of related literature.  In other
contexts, the spherical subgroups are of interest because they are
related to multiplicity free modules (see the discussion below of work
by Kac \cite{kac} and Arzhantsev \cite{arzhantsev}).

\paragraph{Irreducible finite orbit modules.}
The \label{irred_finite_orbit_modules} next family of examples is of a
representation theoretic nature.  Let $X$ be a closed, connected
subgroup of $G = \GL(V)$ with $V$ a finite dimensional vector space,
irreducible under $X$.  Then $X$ is a reductive group and we may write
$X=X'Z$ where $X'$ and $Z$ are the derived subgroup and the center of
$X$ respectively.  Then $Z$ acts as scalars upon $V$.  To make the
following statement easier, we assume that $Z$ equals the full group
of scalars.  Let $P_1$ be the stabilizer in $\GL(V)$ of a $1$-space in
$V$.  We identify $\GL(V)/P_1$ with the collection of $1$-spaces in
$V$.

\begin{theorem}
Let $X$, $\GL(V)$ and $P_1$ be as described.  Then $X$ has a finite
number of orbits upon $V$ if and only if $\double {X'}{\GL(V)}{P_1}$
is finite.  Moreover, all instances where $\double {X'}{\GL(V)}{P_1}$ is
finite have been classified.
\end{theorem}
The classification was found first by Kac \cite{kac} in characteristic
$0$ and extended to positive characteristics by Guralnick, Liebeck,
Macpherson and Seitz \cite{GLMS}.  The latter paper also classifies
when $\double {X'}{\GL(V)}{P_i}$ is finite where $P_i$ is the
stabilizer of an $i$-dimensional subspace.  In Kac's paper, he
suggests a broader context for this problem.  Essentially, he wants to
find those representations which allow a classification of orbits and
he lists various properties which such representations might have in
common.  For instance, one of the properties he studies is that the
modules be spherical (i.e. a Borel subgroup has a dense orbit on $V)$
which is equivalent to the representation being multiplicity free.  A
related problem was studied by Arzhantsev \cite{arzhantsev}.  Here the
goal was to classify those representations of $G$ on $V$ such that a
Borel subgroup has an open orbit in $Gv$ for each $v\in V$.  Again, a
classification is obtained.  This demonstrates Kac's heuristic that
representations with certain nice properties should be amenable to
classification.  The question remains as to which properties are
important.  The present paper takes a small step towards making the
case that having a finite double coset collection is an important
property.

\paragraph{A program for double cosets.}
The last two families of examples share a theme.  In each case we have
reductive groups acting on $G / P$ where $P$ is a parabolic subgroup. 
In an effort to classify all instances where $\double XGY$ is finite
this setting is a natural place to concentrate for reasons we will
discuss.  Firstly, one should start by classifying finiteness of
$\double XGY$ where both $X$ and $Y$ are maximal.  Then the Borel-Tits
Theorem (see \cite[30.4]{humphreys}) implies that each of $X$ and $Y$
is reductive or parabolic.  If $X$ and $Y$ are both parabolics, then
the collection is finite by the Bruhat decomposition.  By work of
Brundan \cite{brundan2}, if $X$ and $Y$ are each reductive and either
maximal or the Levi factor of a parabolic subgroup, then finiteness of
$\double XGY$ is equivalent to having a factorization $G=XY$.  Such
factorizations are rare and have been essentially classified in a
paper by Liebeck, Saxl and Seitz \cite{LSS}.  (Note, having such a
factorization is equivalent to having $\sdouble XGY=1$, so this is
also an example of a finite double coset problem.)  This leaves the
case where $X$ is reductive and $Y=P$ is a parabolic.  (Further
discussion of these matters may be found in the article by Seitz
\cite{seitz}.)

Within the setting of $X$ being reductive and $P$ being parabolic, the
question becomes which reductive groups should be addressed first.  In
many problems, the maximal rank reductive subgroups are the most
important.  Therefore in this paper we will study the case where $X$
is maximal rank.  We intend to classify all cases where $X$ is maximal
rank and reductive and $\double XGP$ is finite in future papers. 
Afterwards, it is hoped that a classification may be obtained where
the only restriction placed upon $X$ will be that it is reductive.

\paragraph{Outline of remaining sections.}
The outline of the rest of this paper is as follows: section
\ref{section_preliminaries} includes basic results and preliminaries;
section \ref{section_reduction_to_finite_groups} reduces the double
coset question of algebraic groups to a related question about finite
groups; section \ref{section_dimension_criterion} applies character
theory to the finite groups (roughly following Lawther \cite{lawther})
and obtains the Dimension Criterion; section
\ref{section_proof_of_sph_or_sph_levi} proves Theorems
\ref{completeness_of_table} and \ref{sph_or_sph_lev}, assuming
Proposition \ref{prop_establishing_complements}; section
\ref{section_constructing_complements} proves Proposition
\ref{prop_establishing_complements}.

\section{Preliminaries}
\label{section_preliminaries}
In this section we list basic results which will be used later.  Many
(perhaps all) of the results in this section are known to others.  We
list them here either for convenience, or because references are
difficult to find.  For standard facts in algebraic groups we refer to
the books by Borel \cite{borel}, Humphreys \cite{humphreys} and
Springer \cite{springer}.

For a reductive, connected reductive group, $G$ we will use the
following notation.  Let $T$ be a maximal torus, $\rtsys$ be the root
system and for each $\alpha\in \rtsys$ let $U_\alpha$ be the
corresponding root group.

\begin{lemma}
Let \label{lemma_A2s_B2s_all_conjugate} $G$ be a simple
algebraic group.  All MRR subgroups of type $A_2$, of the same length,
are conjugate.  All MRR subgroups of type $B_2$ are conjugate.  If the
rank of $G$ is at least three then these subgroups are all Levi
factors of parabolic subgroups.
\end{lemma}

\begin{proof} 
The last statement is clear.  Let $H$ and $H'$ be two MRR subgroups
which are claimed to be conjugate.  By conjugation we may assume that
$H$ and $H'$ share a common maximal torus $T$.  If the rank of $G$ is
two then $H$ and $H'$ are equal.  Otherwise $H$ and $H'$ are Levi
factors and each is generated by $T$ and the root groups (positive and
negative) corresponding to a pair of adjacent nodes in the Dynkin
diagram of $G$.  Then $H$ and $H'$ are conjugate by the action of the
Weyl group.
\end{proof}

The following lemma allows us to make a variety of easy reductions, or 
assumptions about $G$.  

\begin{lemma}
Let \label{lemma_invariants_double_coset} $G$ be a group,
$X$ and $P$ subgroups.  Let $Z$ be the center of $G$, suppose $Z\le P$
and let $\overline X$, $\overline G$ and $\overline P$ be the images
of $X$, $G$ and $P$ under the map $G\to G/Z$.  Let $K$ be a finite
normal subgroup of $G$ and let $\hat X$, $\hat G$ and $\hat P$ be the
images of $X$, $G$ and $P$ under the map $G\to G/K$.  Let $g,h \in G$.
The following are equivalent:
\begin{theorem_count}
\item $\sdouble XGP<\infty$.
\item $\sdouble {\hat X}{\hat G}{\hat P}<\infty$.
\item $\sdouble {\overline X}{\overline G}{\overline P}<\infty$.  
\item $\sdouble {X^g}{G}{P^h}<\infty$.
\end{theorem_count}
\end{lemma}
\begin{proof} 
These statements can all be proven in an elementary fashion.
\end{proof}

Brundan \cite{brundan} also states most of these.  This lemma shows
that the question of whether $\double XGP$ is finite depends only upon
the Lie type of the groups involved.  In particular, it does not
depend upon which elements of an isogeny class are chosen, the
presence of centers, connectedness etc.  By this result, we may work
with either $\GL(V)$ or $\SL(V)$ and get essentially the same results. 
We may also assume that $G$ has simply connected derived subgroup
which eases some of the proofs.  Finally, if $X$ and $P$ are maximal
rank then we may assume that they contain a common maximal torus.

\paragraph{Convention.}
\label{conventions1}
If $\tau$ is an endomorphism we denote by $G_\tau$ the fixed points of
$\tau$ in $G$.  If $G$ is a group and $g \in G$ then $G_g$ denotes the
centralizer of $g$ in $G$.  Finally, $G_{\tau,g}$ denotes those points
in $G$ fixed by both $\tau$ and $g$.  The finite groups of Lie type
arise as the fixed points in $G$ of a Frobenius morphism
$\sigma:G\to G$, where $G$ is defined over the algebraic closure
$\overline\F_p$ of the field $\F_p$ of $p$ elements.  We refer to
\cite{carter} and \cite{steinberg} for details.

\begin{lemma}
Let \label{lemma_dense_subset_center} $G$ be a connected reductive
group with simply connected derived subgroup.  Let $T$ be a maximal
torus of $G$.
\begin{theorem_count}
\item The center of $G$ is contained in each maximal torus of $G$.
(This does not require that $G$ have simply connected derived subgroup.)

\item If $s\in G$ is semisimple then $G_s$ is reductive and connected.

\item For each $s\in T$ we have $Z(G_s)\le T$.

\item The set $\{G_s \st s\in T\}$ is finite.  Its size may be bounded 
by a constant depending only upon the root system of $G$.

\item Fix $s\in T$.  There exist $t_1,\dots,t_r\in T$ such that $\{G_t
\st t\in G,\ G_t > G_s\} = \{G_{t_i} \st 1\le i \le r\}$.  Let
$Z(s)=\{t\in G \st G_t = G_s\}$.  Then $Z(s)$ is an open subset of
$Z(G_s)$ and its complement is $\bigcup_i Z(G_{t_i})$.  

\item If $S$ is a torus and $L=\C G(S)$ then $\{s\in S \st L=G_s\}$ is 
a dense subset of $S$.
\end{theorem_count}
\end{lemma}

\begin{proof}
Part (i) is \cite[26.2]{humphreys}.  

Part (ii) is \cite[3.5.4,3.5.6]{carter}.

Part (iii).  Note that $T$ is a maximal torus of $G_s$.  By part (ii)
we may apply part (i) to the group $G_s$.

Part (iv).  By \cite[3.5.3]{carter} we have that $G_s$ is generated by
$T$, the root groups it contains and by certain elements of the Weyl
group.  Since the Weyl group is finite and the number of root groups 
is finite, there are finitely many possibilities for $G_s$.

Part (v).  Let $U$ be the set proposed as the complement of $Z(s)$. 
Clearly $G_t > G_s \implies Z(G_t) < Z(G_s) \le T$. 
This, and the previous part, show that $t_1,\dots,t_r$ may be chosen
in $T$ as stated and that $U\subseteq Z(G_s)$.  Since the union is
finite, $U$ is a closed set.  Given $t\in Z(G_s)$ we have:
$$\begin{array}{ll}
t\not \in Z(s) & \iff G_t > G_s \\
&\iff G_t = G_{t_i} \text{ for some } i \\
&\iff t\in Z(G_{t_i}) \text{ for some } i\\
&\iff t \in U.
\end{array}$$
This shows that $U$ is the desired complement.

Part (vi).  Note that $S \le Z(L)$.  Then for all $t\in S$ we
have $G_t \ge L$.  By an argument similar to that for part
(v), one can show that the set of $t\in S$ with $G_t > L$ is a
proper, closed subset of $S$.
\end{proof}

\begin{lemma}
Let \label{lemma_centers_rational_centralizers} $G$ be a connected
reductive group and $\sigma :G\to G$ a Frobenius morphism.  Then
$Z(G_\sigma) = Z(G)_\sigma$.  Moreover, if $G$ has simply connected
derived subgroup and $s,t \in G_\sigma$ are semisimple elements, then
$G_s=G_t \iff G_{\sigma,s}= G_{\sigma,t}$.
\end{lemma}

\begin{proof}
Carter \cite[3.6.8]{carter} shows that $Z(G_\sigma)=Z(G)_\sigma$.  For
the second statement note that ``$\Rightarrow$'' is obvious, for
``$\Leftarrow$'' suppose $G_{\sigma,s}= G_{\sigma,t}$.  By Lemma
\ref{lemma_dense_subset_center}(ii) we have that $G_s$ and $G_t$ are
connected and reductive.  Then the first statement shows that
$Z(G_{\sigma,s})= Z(G_s)_\sigma$ whence $t$ is in
$Z(G_{\sigma,s})=Z(G_s)_\sigma$.  This shows that $t$ is in $Z(G_s)$
whence $G_t\ge G_s$.  A symmetric argument shows that
$G_s\ge G_t$.
\end{proof}

\begin{lemma}[Rational normalizer theorem]
Let $G$ be a connected reductive group defined over $\overline \F_p$, 
$\sigma:G\to G$ a Frobenius morphism.  Let $P$ be a $\sigma$-stable parabolic 
subgroup.  Then $N_{G_\sigma}(P_\sigma)= P_\sigma = (N_G(P))_\sigma$
\end{lemma}

\begin{proof}
It is well known that $P=N_G(P)$, which gives the second equality. 
For the first equality it is clear that $P_\sigma \leq N_{G_\sigma}
(P_\sigma)$.  The reverse inclusion follows from the fact that if
$\tilde P$ is a $\sigma$-stable parabolic subgroup with $\tilde
P_\sigma = P_\sigma$ then $\tilde P=P$, see \cite[4.20]{borel_tits}.
%
%
\end{proof}

\begin{cor}
Let \label{rational_points_and_induced_trivial_character}
$G$ be connected and reductive, $P$ a parabolic
subgroup and $x\in G_\sigma$.  Suppose that $\sigma$ is a Frobenius morphism
of $G$ which fixes $P$ and $x$.  Let $(G/P)_x$ be the variety
of $G$-conjugates of $P$ which contain $x$.  Then $\sigma$ acts upon
$(G / P)_x$ and $1_{P_\sigma}^{G_\sigma}(x)$ is equal to the
number of $\sigma$-fixed points on this variety.
\end{cor}

\begin{proof}
Using the Lang-Steinberg Theorem \cite{springer-steinberg} it is easy
to show that $\phi : {G_\sigma}/ {P_\sigma}\to ( G / P)_\sigma$ taking
$gP_\sigma$ to $gP$ is an $x$-equivariant bijection.  Together with
the rational normalizer theorem this shows that we have bijections
between $(G/P)_{\sigma,x}$, $(G_\sigma/P_\sigma)_x$ and
$\{{}^gP_\sigma \st g\in G_\sigma,\ x\in {}^gP_\sigma\}$.  Elementary
character theory shows that $1_{P_\sigma}^{G_\sigma}(x)$ equals the
size of the last collection.

%
%

%

\end{proof}


\begin{lemma}[{\cite[3.5]{nori}}]
Let \label{bounding_rational_points_of_groups}
$H$ be a connected algebraic group of dimension $d$ and 
$\sigma:H\to H$ a standard $q\th$ power Frobenius map.  Then $(q-1)^d 
\le |H_\sigma| \le (q+1)^d$.  
\end{lemma}

\begin{lemma}
Let \label{bounding_rational_points_of_dense_subset_of_Z} $G$ be a
connected reductive group with simply connected derived subgroup,
$\sigma:G\to G$ a standard $q\th$ power Frobenius map and $T$ a
$\sigma$-stable maximal torus.  Fix $s\in T_\sigma$ and let $Z(s)$ and
$t_1$, \dots, $t_r$ be as in Lemma \ref{lemma_dense_subset_center}.  Let
$c_1$ and $d_2$ be the number of components and dimension of $Z(G_s)$
respectively.  Let $I\subseteq \{1,\dots,r\}$ such that $\dim 
Z(G_{t_i}) < \dim Z(G_s) \iff i\in I$.  Let $m=|I|$, let $c_2$ and 
$d_2$ be the maximal number of components and the greatest dimension 
of the $Z(G_{t_i})$ with $i\in I$.  Note that if $m\ne0$ then 
$d_2 < d_1$.  

Then $Z(s)$ is $\sigma$-stable  and
$$(q-1)^{d_1}-mc_2 (q+1)^{d_2}
\ \le \ 
|Z(s)_\sigma| 
\ \le \ 
c_1(q+1)^{d_1}
.$$
\end{lemma}


\begin{proof}
If $t\in Z(s)$ then $G_{\sigma(t)} = \sigma(G_t) = \sigma G_s = G_s$
which shows that $Z(s)$ is $\sigma$-stable.  We will be using Lemma
\ref{bounding_rational_points_of_groups} throughout the proof without 
further mention.

Since $Z(s)\subseteq Z(G_s)$ we see that $|Z(s)_\sigma | \le |
Z(G_s)_\sigma| \le c_1(q+1)^{d_1}$ where the second inequality is
found by calculating $|Z(G_s)_\sigma|$ under the assumption that
$\sigma$ stabilizes each component of $Z(G_s)$.

From Lemma \ref{lemma_dense_subset_center} we have that $Z(G_s) =
Z(s)\ \sqcup\ \bigcup_{i\ge 1} Z(G_{t_i})$, where ``$\sqcup$''
indicates a disjoint union.  Intersecting with $sZ(G_s)^\circ$ and 
taking fixed points we have
$$|(sZ(G_s)^\circ \cap Z(s))_\sigma| 
   = |(sZ(G_s)^\circ)_\sigma| 
   -
   \Big|\Big(sZ(G_s)^\circ \cap (\bigcup_{i\ge 1} 
   Z(G_{t_i}))\Big)_\sigma\Big|.$$
It is easy to check that:
$$\begin{array}{l}
|(sZ(G_s)^\circ \cap Z(s))_\sigma|\ \le \ |Z(s)_\sigma|, \\
(q-1)^{d_1}\ \le \ |Z(G_s)^\circ_\sigma|\ =\ |(sZ(G_s)^\circ)_\sigma|, \\
\Big|\Big(sZ(G_s)^\circ\cap (\bigcup_{i\in I} Z(G_{t_i})) \Big)_\sigma\Big| 
      \ \le\  mc_2(q+1)^{d_2}. 
\end{array}$$
To finish we will prove that $sZ(G_s)^\circ \cap (\bigcup_{i\ge 1}
Z(G_{t_i})) = sZ(G_s)^\circ \cap (\bigcup_{i \in I} Z(G_{t_i}))$.  It
suffices to show that $sZ(G_s)^\circ \cap Z(G_{t_i})$ is empty if
$\dim Z(G_{t_i}) = \dim Z(G_s)$.  Let $\dim Z(G_{t_i})= \dim Z(G_s)$. 
Then $Z(G_{t_i})^\circ = Z(G_s)^\circ$ and $sZ(G_s)^\circ \cap
Z(G_{t_i})$ is empty or all of $sZ(G_s)^\circ$.  However, by
definition of the $t_i$, we have $s\not \in Z(G_{t_i})$ so we are
done.
\end{proof}

%
%
%

\section{Reduction to finite groups}
\label{section_reduction_to_finite_groups}
In this section we prove the fundamental results which relate double
cosets in algebraic groups to double cosets in finite groups.  These
results seem intuitive, but use material surprisingly far from group
theory.  By a reduced algebraic group scheme over $\Z$, we mean
(naively) that the group $G$ is defined, as a subgroup of $\GL_n(\Z)$,
with a finite number of polynomials over $\Z$ and that $\Z[G]$ has no
nilpotents except $0$.  This is the case for the simple algebraic
groups, as well as their parabolic subgroups and generic MRR subgroups
(see  \cite{demazure} or \cite{jantzen}).  Such a group scheme has a
group of points over every field.  For an algebraically closed field
$\k$ one may identify the group of points (of the group scheme) over
$\k$ with the algebraic group (in the naive sense) over $\k$.  The 
field $\overline \F_p$ is the algebraic closure of the field of $p$ 
elements for the prime $p$.  Lawther \cite{lawther} also uses some of 
these results.  

\begin{prop} 
\label{reduction_to_pos_chars}
Let $G$ be a simple algebraic group scheme and $X$ and $P$ closed
algebraic subgroup schemes of $G$, all of which are reduced over $\Z$. 
For a field $\F$ we denote by $G(\F)$, $X(\F)$ and $P(\F)$ the group
of points over $\F$ of $G$, $X$ and $P$.  Let $\k$ be an
algebraically closed field.
\begin{theorem_count}
\item
If $\charistic \k =0$ then 
$$\sdouble {X(\k)}{G(\k)}{P(\k)} <
\infty \iff {\displaystyle \limsup_{p\to \infty}} \sdouble
{X(\overline \F_p)}{G(\overline \F_p)} {P(\overline \F_p)} \
<\infty.$$ 
\item If $\charistic \k=p > 0$ then 
$$\sdouble {X(\k)}{G(\k)}{P(\k)} <
\infty \iff \sdouble {X(\overline \F_p)}{G(\overline
\F_p)}{P(\overline \F_p)}<\infty.
$$
\end{theorem_count}
\end{prop}

\begin{proof}
Part (ii) is proven in \cite{GLMS}.  (We view the group $X(\k)\times
P(\k)$ as acting on the affine space $G(\k)$.  The assumption in
\cite{GLMS} that $X(\k)\times P(\k)$ should be reductive is not used.) 
It may also be proven using a model theoretic argument similar in
nature to the one we give now for part (i).  For basic facts about
model theory we refer to the textbooks by Fried-Jarden \cite{FJ} or
Hodges \cite{hodges}.

For $p$ equal to $0$ or a prime, let $ACF_p$ be the theory of
algebraically closed fields of characteristic $p$.  Then $ACF_p$ is a 
complete theory. 

For a field $\F$ we identify $G(\F)$ as a set of matrices in
$\GL_n(\F)$ using the defining polynomials over $\Z$.  We make similar
identifications for $X$ and $P$.  Since $G$, $X$ and $P$ are defined
over $\Z$ we can express membership in $G(\F)$, $X(\F)$ and $P(\F)$
with first order sentences.  Let $\phi$ be the sentence which, applied
to the model $\F$, gives $\exists g_1$, \dots, $g_n\in G(\F)$,
$\forall g\in G(\F)$, $\exists x\in X(\F)$, $\exists y\in P(\F)$,
$\exists i\in \{ 1,\dots, n \}$ such that $xgy=g_i$.  In other words,
$\phi$ applied to $\F$ states that $\sdouble {X(\F)}{G(\F)}{P(\F)}\le
n$.

Suppose $\double {X(\k)} {G(\k)} {P(\k)}$ is infinite in
characteristic zero.  Then $\phi$ is false in $\k$.  Then $ACF_0\vdash
\neg \phi$ by completeness.  This means that we have a finite number
of steps each of which only uses a finite number of axioms to derive
$\neg \phi$.  In particular, only finitely many axioms which assert
that $m\cdot 1\ne 0$ are used and so there exists a prime $p_0$ which
is greater than every $m$ which is used in this manner. 
For all primes $p\ge p_0$ the axioms and steps which are used in the
proof of $ACF_0\vdash \neg \phi$ may also be used to conclude
$ACF_p\vdash \neg \phi$.  Therefore, for all such $p$ we have
$\sdouble{X(\overline \F_p)} {G(\overline\F_p)} {P(\overline \F_p)}>
n$ whence $\limsup_{p\to \infty}\sdouble{X(\overline\F_p)}
{G(\overline\F_p)} {P(\overline\F_p)}> n$.

Conversely, a similar argument shows that $ACF_0 \vdash \phi \implies
ACF_p\vdash \phi$ for all $p$ sufficiently large.  Therefore
finiteness in characteristic $0$ implies boundedness of
$\sdouble{X(\overline\F_p)}{G(\overline\F_p)}{P(\overline\F_p)}$ as
$p\to \infty$.
\end{proof}

\begin{lemma} Let \label{exceptional_main_reduction_theorem}
$G$ be a connected algebraic group defined over $\k= \overline \F_p$,
$\sigma:G\to G$ a Frobenius morphism, $X$ and $P$ closed subgroups 
which are $\sigma$-stable.  If $\double XGP$ is infinite let $C=1$.  If 
$\double XGP$ is finite let $C$ be an upper bound on the number of 
components of stabilizers of $X\times P$ acting on $G$.  Then
$$
\frac 1C  \limsup_{n\to \infty}
\sdouble {X_{\sigma^n}} {G_{\sigma^n}} {P_{\sigma^n}} 
\quad \le \quad 
\sdouble XGP
\quad \le \quad 
 \limsup_{n\to \infty} \sdouble  {X_{\sigma^n}} {G_{\sigma^n}} {P_{\sigma^n}}.
$$
\end{lemma}

\begin{proof}
Suppose ${\displaystyle \limsup_{n\to \infty} \sdouble {X_{\sigma^n}}
{G_{\sigma^n}} {P_{\sigma^n}}}$ is finite and less than $m$.  We will
show that $\sdouble XGP<m$.  Let $g_1,\dots,g_m \in G$.  There is a
number $n\in \N$ such that $g_1,\dots,g_m \in G_{\sigma^n}$ and
$m>\sdouble {X_{\sigma^n}} {G_{\sigma^n}} {P_{\sigma^n}}$.  Then at
least two of $g_1,\dots,g_m$ are in the same $X_{\sigma^n}\times
P_{\sigma^n}$-orbit, whence they are in the same $X\times P$-orbit. 
Since this holds for every $g_1,\dots,g_m \in G$ we see that $\sdouble
XGP<m$.

Suppose now that $\double XGP$ is finite, let $n$ be given and let
$(\double XGP)_{\sigma^n}$ be the collection of $\sigma^n$-stable
$X\times P$-orbits.  Then the Lang-Steinberg Theorem
\cite{springer-steinberg} shows that $C\sdouble XGP \ge
C\size{(\double XGP)_{\sigma^n}} \ge \sdouble {X_{\sigma^n}}
{G_{\sigma^n}} {P_{\sigma^n}}$.
\end{proof}

%

\section{Character Theory and the Dimension Criterion}
\label{section_dimension_criterion}
\paragraph{Strategy and conventions.}
By Lemma \ref{lemma_invariants_double_coset} we may assume that $G$
has simply connected derived subgroup when convenient and use Lemma
\ref{lemma_dense_subset_center}, Lemma
\ref{lemma_centers_rational_centralizers} etc.  We wish to establish a
criterion for infiniteness which is independent of characteristic.  By
section \ref{section_reduction_to_finite_groups} it suffices to work
over the algebraic closures of finite fields and establish
infiniteness independently of the field.  Let $G$ be defined over the
algebraic closure of a field of positive characteristic.  Let
$\sigma:G\to G$ be a $q\th$ power Frobenius morphism.  We assume that
$X$ and $P$ are $\sigma$-stable.  Then to prove infiniteness it
suffices to show that $\sdouble {X_{\sigma^n}} {G_{\sigma^n}}
{P_{\sigma^n}}$ is unbounded as $n$ approaches infinity.  For fixed
points we will use the notation $G_\sigma$, $P_s$, etc as described in
section \ref{conventions1}.  If $H$ is a group, the notation
$[h]\subseteq H$ means $h$ is an element of $H$ and $[h]$ is its
$H$-conjugacy class.  An element denoted by $s$ will be semisimple,
and an element denoted by $u$ will be unipotent.  This development
roughly follows Lawther \cite{lawther}.  Essentially the following
lemma regroups the terms in the inner product of characters.

\begin{lemma}
Let \label{lemma_grouping_terms_in_inner_product} $X$ and $P$ be
subgroups of $G$ with $P$ parabolic.  Define
\label{first_reduction_to_charactcers} an equivalence relation on
semisimple elements in $X_\sigma$ as follows: $s$ and $t$ are
equivalent if $G_{\sigma,s}$ and $G_{\sigma,t}$ are
$X_\sigma$-conjugate.  Denote the equivalence class of $s$ by
$E(s,\sigma)$.  Choose a set $S_\sigma$ of representatives of these
equivalence classes.  Then
$$
\sdouble{X_\sigma}{G_\sigma}{P_\sigma}
= 
\sum_{s\in S_\sigma}
\sum_{[u]\subseteq X_{\sigma,s}}
\sizefrac{E(s,\sigma)}{X_\sigma} 
\sizefrac{X_{\sigma,s}}{X_{\sigma,s,u}}
 1_{P_\sigma}^{G_\sigma}(su).
$$
\end{lemma}
\begin{proof}
Basic character theory gives 
$$
\sdouble {X_\sigma}{G_\sigma}{P_\sigma}
    =  (1_{X_\sigma}^{G_\sigma}, 1_{P_\sigma}^{G_\sigma})_{G_\sigma}^{} 
    =    (1_{X_\sigma}, 1_{P_\sigma}^{G_\sigma})_{X_\sigma}^{} 
    =  \frac 1{\size {X_\sigma}}\sum_{x\in X_\sigma}1_{P_\sigma}^{G_\sigma}(x).
$$
Applying the Jordan-Chevalley decomposition within the finite group
$X_\sigma$ we get that this last sum is equal to
$$\frac 1{\size {X_\sigma}}
    \sum_{s\in X_\sigma}\,
    \sum_{ u\in X_{\sigma,s}}  1_{P_\sigma}^{G_\sigma}(su).$$ 
Now we claim that $t\in E(s,\sigma)$ implies that
$$
\sum_{ u\in X_{\sigma,t}}
1_{P_\sigma}^{G_\sigma}(tu)
=
\sum_{ u\in X_{\sigma,s}}
1_{P_\sigma}^{G_\sigma}(su).$$ 
Let $x\in X_\sigma$ with $(G_{\sigma,t})^x = G_{\sigma,s}$.  The
crucial step is to show that for all $u\in X_{\sigma,t}$ we have
$1_{P_\sigma}^{G_\sigma}(tu) = 1_{P_\sigma}^{G_\sigma}(su^x)$.  Once
this is done, conjugation by $x$ shows that the sums are equal.  We
work at the level of algebraic groups.  Given $u\in X_{\sigma,t}$, let
$( G/P)_{ tu}$ and $( G/P)_{ su^x}$ be the varieties of conjugates of
$P$ which contain $tu$ and $su^x$ respectively.  Then, by Lemma
\ref{rational_points_and_induced_trivial_character}
$1_{P_\sigma}^{G_\sigma}(tu)$ and $1_{P_\sigma}^{G_\sigma}(su^x)$ are
the numbers of $\sigma$-rational points on these varieties.  Let $g\in
G$ such that $t\in P^g$, let $T$ be a maximal torus of $P^g$ which
contains $t$.  We apply Lemma
\ref{lemma_centers_rational_centralizers} to see that $(G_t)^x = G_s$.  We have
the following:
$$T\leq G_t \implies T^x \leq G_s
\implies s\in T^x \implies
s\in P^{gx}.$$
It is now easy to see that $tu\in
P^g\implies su^x\in P^{gx}$.  Therefore, conjugation by $x$ gives a
$\sigma$-equivariant bijection $(G/P)_{ tu}\to (
G/P)_{ su^x}$.  Taking $\sigma$-fixed points and applying 
Lemma \ref{rational_points_and_induced_trivial_character} finishes the claim.

Using the claim we have 
$$  \frac 1{\size {X_\sigma}} 
    \sum_{ s\in X_\sigma}\,
    \sum_{ u\in X_{\sigma,s}}  
    1_{P_\sigma}^{G_\sigma}(su) 
    =
    \frac 1{\size {X_\sigma}}
    \sum_{ s\in S_\sigma}
    \sum_{ u\in X_{\sigma,s}}
    |E(s,\sigma)|
    1_{P_\sigma}^{G_\sigma}(su).$$ 
To finish the proof note that if $u$ and $u'$ are conjugate in
$X_{\sigma,s}$ then $1_{P_\sigma}^{G_\sigma}(su)=
1_{P_\sigma}^{G_\sigma}(su')$.  Therefore the sum over unipotent
elements in $X_{\sigma,s}$ may be replaced by representatives of
unipotent classes.  The size of such a unipotent class is
$$\dsize { u^{ X_{\sigma,s}}}=
\sizefrac{X_{\sigma,s}}{X_{\sigma,s,u}}.$$
\end{proof}

\begin{lemma}
Let $X$ and $P$ be subgroups of $G$ with $X$ maximal rank.
\label{calc_size_eq_class}
Let $T$ be a $\sigma$-stable maximal torus in $X$, $W=N_G(T)$ the Weyl
group, and $s\in T_\sigma$.  Let $Z(s,\sigma)= \{ t \in G_{\sigma} \st
G_{\sigma,t} = G_{\sigma,s}\}$, $Z(s)=\{ t \in G\st G_t = G_s\}$ (as
in Lemma \ref{lemma_dense_subset_center}) and let $E(s,\sigma)$ be as in the
previous lemma.  Then we have:
\begin{theorem_count}
\item $Z(s,\sigma)= Z(s)_\sigma$.\medskip
\item $\displaystyle
\begin{matrix}\frac1 {\size W} \size{s^{X_\sigma}\times Z(s,\sigma)} & \le & 
\size{E(s,\sigma)}
& \le & \size{s^{X_\sigma}\times Z(s,\sigma)}.
\end{matrix}$
\end{theorem_count}
\end{lemma}
\begin{proof}
Part (i).  If $t\in Z(s)_\sigma$ then $G_s=G_t$ which implies
$G_{\sigma,s}= G_{\sigma,t}$, whence $Z(s)_\sigma \subseteq
Z(s,\sigma)$.  If $t\in Z(s,\sigma)$ then $G_{\sigma,s}=
G_{\sigma,t}$.  By Lemma \ref{lemma_centers_rational_centralizers}  
this implies $G_s = G_t$ whence $t\in Z(s)$.

Part (ii).  It suffices to show that the following is a surjective map, with
fibers bounded in size by $\size W$: $$\phi:\map{s^{X_\sigma}_{}\times
Z(s,\sigma)} {E(s,\sigma)} {[2\jot](s^x,t)} {t^x}.$$ Note that this
map is well-defined as every element in $X_\sigma$ which centralizes $s$
also centralizes $t$.  To see that the map is surjective, let $t\in
E(s,\sigma)$ and let $ x\in X_\sigma$ with $(G_{\sigma,t})^x=
G_{\sigma,s}$.  Then $(s^{x\inv},t^x)$ is in the domain of $\phi$ and
$\phi(s^{x\inv },t^x)=t$.

The remainder of the proof bounds the size of the fiber.  Let
$(s^x,t_1)$ be an element of the domain.  We claim that
$$\phi\inv (t_1^x) 
= 
\left\{(s^{w\inv x},t_1^w) \st w\in W \right \}\ \bigcap\
\Bigl ( s^{X_\sigma}\times Z(s,\sigma) \Bigr ).$$ 
It is easy to see that the set on the right is contained in
$\phi\inv(t_1^x)$.  Conversely, let $(s^y,t_2)\in \phi\inv(t_1^x)$ and
note $t_2= t_1^{xy\inv}$.  

Claim: $T$ contains $t_1$, $t_2=t_1^{xy\inv}$, $s$, and $s^{yx\inv}$. 
By assumption $s\in T$.  Note that Lemmas
\ref{lemma_dense_subset_center} and
\ref{lemma_centers_rational_centralizers} show $Z(G_{\sigma,s})\leq
T$.  By definition $t_1,t_2 \in Z(G_{\sigma,s})$.  The following
calculation shows that $s^{yx\inv}\in Z(G_{\sigma,s})$ (which finishes
the claim):
$$ 
G_{ \sigma,s^{yx\inv}} = (G_{\sigma,s})^{yx\inv}=
(G_{\sigma,t_2})^{yx\inv} = G_{ \sigma,t_2^{yx\inv}} = G_{\sigma,t_1} =
G_{\sigma,s}.$$ 

We will use \cite[3.7.1]{carter}, and the (standard) notation which
appears there.  Write $xy\inv =ut\dot wu'$ in the Bruhat canonical
form.  Since $t_2=t_1^{xy\inv}$ we have $t_2=t_1^w$.  Since $xy\inv $
conjugates $s^{yx\inv}$ to $s$ we have $(s^{yx\inv })^w=s$ and
$s^y=s^{w\inv x}$.  Therefore $(s^y,t_2)=(s^{w\inv x}, t_1^w)$.
\end{proof}


\begin{cor}
Let \label{corollary_sum_with_grouping_and_size_of_eq_class}
$X$ be a maximal rank subgroup of $G$ and $P$ a parabolic subgroup
of $G$.  Let $S_\sigma$ and $Z(s,\sigma)$ be as in
Lemmas \ref{lemma_grouping_terms_in_inner_product} and
\ref{calc_size_eq_class} respectively.  We have $$\frac 1{\size{W}}
\sum \sizefrac {Z({s,\sigma})} {X_{\sigma,s,u}} 1_ {P_{\sigma}} ^
{G_{\sigma}}(su) \quad \le \quad
\sdouble{X_\sigma}{G_\sigma}{P_\sigma} \quad \le \quad \sum \sizefrac
{Z({s,\sigma})} {X_{\sigma,s,u}}
1_ {P_{\sigma}} ^ {G_{\sigma}}(su),$$
where each  sum is taken over $s\in S_\sigma$, $[u]\subseteq X_{\sigma,s}$.
\end{cor}

\begin{proof}
Combine Lemma \ref{first_reduction_to_charactcers} and the bounds for
$E(s,\sigma)$ just obtained in Lemma \ref{calc_size_eq_class}.  Note that
$$\sizefrac {s^{X_\sigma}\times Z(s,\sigma)}{X_\sigma} \sizefrac
{X_{\sigma,s}}{X_{\sigma,s,u}} =
\sizefrac {Z({s,\sigma})} {X_{\sigma,s,u}}.$$ 
\end{proof}

\begin{proofof}[Proof of the Dimension Criterion]
We \label{proof_of_dimension_criterion} have that $X$ is maximal rank,
$P$ parabolic and $s\in X\cap P$.  We assume $\dim Z(G_s) + \dim G_s >
\dim X_s + \dim P_s$.  We also assume that $X$ and $P$ are
$\sigma$ stable where $\sigma$ is a standard $q\th$ power Frobenius
map.

Let $\overline{G_s}$, $\overline{X_s}$ and $\overline{P_s}$ denote the
quotients of $G_s$, $X_s$ and $P_s$ by $Z(G_s)$.  The dimension
inequality we have assumed implies that $\dim G_s - \dim Z(G_s) > \dim
X_s - \dim Z(G_s) + \dim P_s - \dim Z(G_s)$.  This implies $\sdouble
{\overline{X_s}} {\overline{G_s}} {\overline{P_s}}=\infty$ whence 
$\sdouble {X_s}{G_s}{P_s}=\infty$ by
Lemma \ref{lemma_invariants_double_coset}.

It remains to show that $\sdouble XGP=\infty$.  Using Corollary
\ref{corollary_sum_with_grouping_and_size_of_eq_class} and Theorem
\ref{exceptional_main_reduction_theorem}, it suffices to show that the
term in Corollary \ref{corollary_sum_with_grouping_and_size_of_eq_class}
corresponding to $s\in S_\sigma$, $1=[u]\subseteq X_{\sigma,s}$ is
unbounded as we replace $\sigma$ with $\sigma^n$ and let $n\to
\infty$.  This term is
$$ \frac 1{\size W} \sizefrac{Z(s,\sigma^n)}{X_{\sigma^n,s}}
1_{P_{\sigma^n}}^{G_{\sigma^n}}(s).$$ 
It is easy to show that $\displaystyle
1_{P_{\sigma^n}}^{G_{\sigma^n}}(s)\ge \sizefrac
{G_{\sigma^n,s}}{P_{\sigma^n,s}}$ whence this term is bounded below by
$$
\frac 1{\size W}
\sizefrac{Z(s,\sigma^n)} {X_{\sigma^n,s}}
 \sizefrac{G_{\sigma^n,s}} {P_{\sigma^n,s}}.$$
Therefore it suffices to show that
$$\limsup_{n\to \infty}\frac 1{\size W}
\sizefrac{Z(s,\sigma^n)}{X_{\sigma^n,s}} \sizefrac{G_{\sigma^n,s}}
{P_{\sigma^n,s}}= \infty.$$
By Lemmas \ref{bounding_rational_points_of_groups} and
\ref{bounding_rational_points_of_dense_subset_of_Z} we have
$$\limsup_{n\to \infty}
\frac 1{\size W}
\sizefrac{Z(s,\sigma^n)}{X_{\sigma^n,s}}
\sizefrac{G_{\sigma^n,s}}{P_{\sigma^n,s}} 
=
\lim_{n\to \infty} C
\mbox{\large $\displaystyle \frac {(q^n)^{\dim Z(G_s)+\dim G_s}}
      {(q^n)^{\dim X_s + \dim P_s}}$},$$ 
for some constant $C$.  It is now easy to see that this limit is
infinite.
\end{proofof}

\section{Proof of Theorems \ref{completeness_of_table} and
\ref{sph_or_sph_lev}}
\label{section_proof_of_sph_or_sph_levi}

Throughout this section $G$ is a simple algebraic group, $X$ a generic
MRR subgroup, $P$ a parabolic subgroup with Levi factor $L$.  Starting
with Proposition \ref{prop_establishing_complements} we will use $H$
for arguments which apply to both $X$ and $L$.

\begin{lemma}
Let \label{calculation_of_dims} $s\in X\cap L$.  If either of the
following holds then $\double XGP$ is infinite:
\begin{theorem_count}
\item $G_s$ is of type $A_2$ and $X_s$ and $L_s$ are tori.

\item $G_s$ is of type $B_2$, $X_s$ is a torus and $L_s$ is of type
$A_1$ or a torus.
\end{theorem_count}
\end{lemma}

\begin{proof}
Using the dimension criterion it suffices to show that
$$\dim Z(G_s) + \frac 12 \dim G_s -\dim X_s - \frac 12 \dim L_s > 0.$$
It is easy to check in each case that the quantity on the left is at least $1$.
\end{proof}

\begin{cor}  
\label{A2_B2_complments_imply_infinite} 
If either of the following hold then $\double XGP$ is infinite:
\begin{theorem_count}
\item $X$ and $L$ have conjugate $A_2$ complements.

\item $X$ has a $B_2$ complement $K$, and for some conjugate $\tilde
K=K^g$ we have that $\tilde K\cap L$ is a MRR subgroup which is a
torus or of type $A_1$.

\end{theorem_count}
\end{cor}

\begin{proof}
If $G$ has rank 2 and (i) or (ii) holds then it is easy to show that
$\double XGP$ is infinite by dimension.  

Assume now that the rank of $G$ is at least $3$.  If (i) holds let $K$
be the $A_2$ complement of $X$.  Using Lemma
\ref{lemma_invariants_double_coset} we replace $P$, if necessary, by a
conjugate so that in (i) the $A_2$ complements of $X$ and $L$
coincide, or so that in (ii), we may take $\tilde K= K$.  Since the
rank of $G$ is at least $3$, we have that $K$ is a Levi factor of a
parabolic (see Lemma \ref{lemma_A2s_B2s_all_conjugate}), whence is of
the form $\C G(S)$ for some torus $S$.  Apply Lemma
\ref{lemma_dense_subset_center}, to see that there exists $s\in S$
with the centralizer of $s$ in $G$ equal to $K$.  We are done by the
previous lemma.
\end{proof}

\begin{cor}
\label{corollary_complements_not_spherical}
If $X$ has an $A_2$ or $B_2$ complement then $X$ is not spherical.
\end{cor}

\begin{proof}
Apply Lemma \ref{A2_B2_complments_imply_infinite}, noting that the
Levi factor for a Borel subgroup is a torus, which has every type of
complement possible.
\end{proof}

\begin{prop}  
\label{prop_establishing_complements}
Let $H$ be a generic MRR subgroup of $G$ which does not appear in
table \ref{spherical_subgroups}.  The following hold, and, in
particular, $H$ has an $A_2$ or $B_2$ complement in all cases.
\begin{theorem_count}
\item 
\label{part_single_root_length_not_in_table_implies_A2_complement}
If $G$ has single root length, then $H$ has an $A_2$ complement.

\item If $H$ is the Levi factor of a parabolic with non-abelian 
unipotent radical then $H$ has an $A_2$, $B_2$ or $G_2$ complement.  

\item Let $G$ equal $B_n$ or $C_n$.

  \begin{theorem_count}
  \item If $n=2$ then $H$ has a $B_2$ complement.  
  
  \item If $G=B_n$ and $H$ equals $D_{n-1}T_1$ or $D_{n_1}D_{n_2}$
  then $H$ has a $B_2$ complement.
  
  \item Otherwise $H$ has an $A_2$ complement.
  \end{theorem_count}

\item Let $G=B_n$.  If $H$ is a Levi factor then there exists a MRR 
subgroup $K$ of type $B_2$ with $H\cap K$ a MRR  subgroup which is
either a torus or of type $A_1$.

\item If $G=F_4$ the maximal possibilities for $H$ are $C_3T_1$,
$A_2\tilde A_2$, $B_3T_1$, $A_1A_1B_2$, $\tilde A_1A_3$, $D_4$.  The
first possibility has a long $A_2$ complement, the next has both long
and short $A_2$ complements, and the rest have short $A_2$
complements.  In particular, if $L$ is a Levi factor for a parabolic
subgroup which is not an end node parabolic, then $L$ has both long
and short $A_2$ complements.
\end{theorem_count}
\end{prop}

The proof of this proposition is delayed until the next section.

\begin{proofof}[Proof of Theorem \ref{completeness_of_table}]
The work of Brundan \cite{brundan} and Theorem \ref{lawthers_theorem}
show that (ii) $\implies$ (i).  Corollary
\ref{corollary_complements_not_spherical} shows that (i) $\implies$
(iii).  Proposition \ref{prop_establishing_complements} shows that
(iii) $\implies$ (ii).
\end{proofof}

\begin{proofof}[Proof of Theorem \ref{sph_or_sph_lev}]
We assume that $X$ and $L$ are not spherical we will show that
$\double XGP$ is infinite.

If $G=G_2$ then by dimension one finds that $X$ non-spherical implies
$\double XGP$ is infinite.  For the remainder of the proof assume
$G\ne G_2$.

Recall our convention that $D_1$ is a $1$-dimensional torus.  If
$(G,X)\neq (B_n, D_{n_1}D_{n_2})$ then let $H_X$ be an $A_2$
complement for $X$ and let $H_L$ be an $A_2$ complement of $L$, of the
same length as $H_X$ (length is only an issue for $F_4$).  If $G=B_2$,
or $(G,X)= (B_n, D_{n_1}D_{n_2})$ then let $H_X$ be a $B_2$ complement
for $X$ and let $H_L$ be a MRR subgroup of type $B_2$ with $L\cap H_L$
a MRR subgroup of type $A_1$ or a maximal torus.  Apply Lemma
\ref{lemma_A2s_B2s_all_conjugate} to see that $H_X$ and $H_L$ are
conjugate.  Apply Lemma \ref{A2_B2_complments_imply_infinite} to see
that $\double XGP$ is infinite.
\end{proofof}

We offer some insight into this result.  Suppose the rank of $G$ is at
least 3 and that $(G,X)\neq (B_n, D_{n_1}D_{n_2})$.  If $X$ and $L$
are both not spherical, then both have $A_2$ complements.  Assume
these $A_2$ complements are conjugate.  Let $T_2$ and $B$ be a maximal
torus and Borel subgroup of $A_2$ respectively.  Then
$\double{T_2}{A_2}{B}$ is infinite, and, in a sense, it may be
``embedded'' in $\double XGP$.  The embedding comes from the expansion
of terms in the inner product of finite group characters corresponding
to $\double XGP$.  Take $s$ such that its centralizer (modulo a
central torus) gives $T_2$, $A_2$ and $B$ in place of $X$, $G$ and
$P$.  The terms which arise by applying Corollary
\ref{corollary_sum_with_grouping_and_size_of_eq_class} to $X$, $G$ and
$P$ and $s$ are the same terms one would obtain by applying Corollary
\ref{corollary_sum_with_grouping_and_size_of_eq_class} to the groups
$\tilde X = T_2$, $\tilde G= A_2$ and $\tilde P=B$.

\section{Proof of Proposition \ref{prop_establishing_complements}}
\label{section_constructing_complements}
We prove parts (i) and (ii) immediately.  Parts (iii)-(v) follow after
Corollary \ref{corollary_sub-maximal_implies_non_spherical}.

Proposition \ref{prop_establishing_complements}(i).  Recall $G$ has single
root length and $H$ is a MRR subgroup which fails to appear in table
\ref{spherical_subgroups}.  Then, by Theorem \ref{lawthers_theorem},
$H$ is not anti-open, whence there exist $\alpha,\beta,\alpha+\beta\in
\rtsys(G) -\rtsys(H)$.  Let $T$ be the maximal torus used to define
these roots and denote the corresponding root groups by $U_\alpha$,
$U_\beta$ etc.  Let $K$ be the group generated by $T$ and $U_{\pm
\alpha}$, $U_{\pm \beta}$.  Then $K$ is an $A_2$ complement for $G$.

Proposition \ref{prop_establishing_complements}(ii).  Recall that $H$
is the Levi factor of a parabolic with non-abelian unipotent radical. 
Let $Q$ be the unipotent radical of the parabolic associated with $H$. 
Fix a maximal torus $T$ in $H$ and let $U_{\alpha}$ and $U_\beta$ be
root groups contained in $Q$ which do not commute.  Let $K$ be the
group generated by $T$ and $U_{\pm \alpha}$, $U_{\pm \beta}$.  Then
$K$ is of type $A_2$, $B_2$ or $G_2$ and it remains to show that
$K\cap H=T$.  Let $Q^-$ be the opposite unipotent radical of $Q$. 
Note that $\rtsys(K) \leq \rtsys(Q) \cup \rtsys(Q^-)$
whence $\rtsys(K)\cap \rtsys(H)=\emptyset$ whence $K\cap H=T$.  Thus 
$K$ is an $A_2$, $B_2$ or $G_2$ complement for $H$.

\begin{lemma}
Let $\phi$ be an irreducible root system in a Euclidean space $\E$
with inner product $(\ ,\ )$.  Let $\hat\phi$ be a proper, closed
subsystem of $\phi$.  Then $\phi-\hat\phi$ spans $\E$ and for each 
$\beta\in \hat\phi$ there exists $\alpha\in \phi-\hat \phi$ with 
$(\alpha,\beta)\ne 0$.
\end{lemma}

\begin{proof}
  Let $n$ be the dimension of $\E$ and fix a Dynkin diagram $\rtbase$
  of $\phi$.  Given $\alpha$, $\beta \in \rtbase$ the path connecting
  $\alpha$ to $\beta$ is the shortest such path and includes $\alpha$
  and $\beta$.  The sum over this path means the sum of each element
  of $\rtbase$ which is contained in the path.  It is easy to check
  that such a sum is itself a root.
  
  For the first conclusion it suffices to show that we have $n$
  independent vectors in $\phi - \hat\phi$.  Since $\hat\phi$ is a
  proper, closed subsystem we have that $\rtbase - \hat\phi$ is
  non-empty.  For each element $\alpha\in\rtbase $ let $\gamma_\alpha$
  be the path connecting $\alpha$ to some element of $\rtbase -
  \hat\phi$.
%
%
We re-index these paths so that for $i\in\{1,\dots,n\}$ the path
$\gamma_i$ contains a node which does not appear in $\gamma_1,\dots,
\gamma_{i-1}$.  Let $\beta_1,\dots,\beta_n$ be the sums over the paths
just constructed.  By the manner in which the paths $\gamma_i$ were
indexed, it is easy to see that $\beta_1,\dots,\beta_n$ are linearly
independent.  By the manner in which the paths were chosen, we may
write each $\beta_i$ as the sum of a root in $\hat\phi$ and a root
outside of $\hat\phi$.  This shows that $\beta_i$ is not in
$\hat\phi$.

For the final conclusion note that $\beta$ is not orthogonal to $\E$, 
whence it is not orthogonal to $\phi-\hat\phi$.
\end{proof}

\begin{cor}
Let \label{corollary_root_system_complements_span} $G$ be a reductive
algebraic group, $H$ a MRR subgroup.  Let $\phi\leq \rtsys(G)$ be
an irreducible root system and let $\phi(H)= \phi \cap \rtsys(H)$. 
Suppose $\hat \phi$ is a closed subsystem of $\phi$ with $\phi >
\hat \phi > \phi(H)$.

\begin{theorem_count}

\item If $\phi$ has single root length then $H$ has an $A_2$
complement (whose length is the same as $\phi$).

\item If $\phi$ is closed, $G=B_n$ and $\hat\phi-\phi(H)$ contains a
short root then $H$ has a $B_2$ complement.

\end{theorem_count}
\end{cor}


\begin{proof}
Fix $\beta \in \hat \phi- \phi(H)$ and assume $\beta$ is short if (ii)
holds.  By the previous lemma there exists $\alpha\in \phi-\hat \phi$
with $(\alpha,\beta)\ne 0$.  Note that $\alpha\ne \pm\beta$ and that
$i\alpha+j\beta\in \rtsys(G) \implies i\alpha+j\beta\in \phi$.  If
$(\alpha,\beta) > 0$ we may replace one root with its negative and
assume that $(\alpha,\beta)<0$, whence $\alpha+\beta\in\phi$.
Since $\alpha \not \in \hat \phi$ and $\beta \in \hat \phi$ we see
that $\alpha + \beta \not \in \hat \phi$.  Similarly, we see that
$\alpha + 2 \beta \not \in \hat \phi$ (of course it may not even be a
root) and that $2 \alpha + \beta$ is not a root.  Let $T$ be a maximal
torus in $H$ used to define these roots.  Set $K$ to be the group
generated by $T$ and all root groups $U_{i\alpha+j\beta}$ where
$i\alpha+j\beta$ is a root with $i$ and $j$ integers.  Then $K\cap
H=T$ and $K$ is of type $A_2$, $B_2$ or $G_2$.  If (i) holds then $K$
is a group of type $A_2$.  If (ii) holds then $K$ is a group of type
$B_2$ since $\beta$ is short and $B_n$ has no MRR subgroups of type
short $A_2$ or $G_2$.
\end{proof}

\begin{cor}  
Let 
\label{corollary_sub-maximal_implies_non_spherical}
$G$ be a reductive algebraic group, $H$ a MRR subgroup, $\phi\leq
\rtsys(G)$ an irreducible root system with single root length.  If
$\phi(H)$ is submaximal in $\phi$ then there exists $\hat \phi$ as in
the previous corollary.
\end{cor}

\begin{proof}
In a root system with single root length, every root subsystem is 
closed.  
\end{proof}

\begin{proofof}[Proof of Proposition \ref{prop_establishing_complements}
(iii),(iv),(v)]

Part (iii).
Recall $G$ equals $B_n$ or $C_n$ and $H$ is a generic MRR subgroup 
which does not appear in table \ref{spherical_subgroups}.  

Part (a).  Since $n=2$ and $H$ does not appear in table
\ref{spherical_subgroups} we see that $H$ is just a torus and $G$
itself is a $B_2$ complement.

Part (b).  We have $n\ge 3$, $G=B_n$ and $H\in\{D_{n-1}T_1,\
D_{n_1}D_{n_2}\}$.  If $H=D_{n_1}D_{n_2}$ then $n_1\ge 2$.  Let
$\phi=\rtsys(G)$ so $\phi(H)= \rtsys(H)$.  If $H=D_{n-1}T_1$ let
$\hat\phi=\rtsys(B_{n-1}T_1)$.  If $H=D_{n_1}D_{n_2}$ let
$\hat\phi=\rtsys(B_{n_1}D_{n_2})$.  In each case $\hat\phi$ contains a
short root and $\phi(H)$ does not.  Then $\hat\phi-\phi(H)$ contains a
short root and we are done by part Corollary
\ref{corollary_root_system_complements_span}(ii).

Part (c).  Note that $n\ge 3$ and if $G=B_n$ that $H\not\in\{
D_{n-1}T_1$, $D_{n_1}D_{n_2}\}$.  If $G=B_n$ let $\phi$ and $\phi(H)$
equal the long roots in $\rtsys(G)$ and $\rtsys(H)$ respectively. 
If $G=C_n$ let $\phi$ and $\phi(H)$ equal the short roots in 
$\rtsys(G)$ and $\rtsys(H)$ respectively.
%
%
In both cases $\phi$ is of type $D_n$, and is irreducible since $n\ge
3$.  The maximal subsystems of $\phi$ are $A_{n-1}$, $D_{n-1}$ and
$D_{n_1}D_{n_2}$.  The subsystem $\phi(H)$ cannot equal $D_n$,
$A_{n-1}$, $D_{n-1}$ or $D_{n_1}D_{n_2}$ as this would contradict
either the assumption that $H$ is not in table
\ref{spherical_subgroups} or the extra restrictions on $H$ when
$G=B_n$.  Therefore $\phi(H)$ is a submaximal subsystem of $\phi$ and
we are done by Corollary
\ref{corollary_sub-maximal_implies_non_spherical}.\medskip

Part (iv).  We have that $G=B_n$ and that $H$ is a Levi factor of $G$. 
Let $T$ be a maximal torus in $H$, let $\alpha_1,\dots,\alpha_n$ be
the nodes in the Dynkin diagram of $G$ in the usual order and suppose
that $H$ is described by ``crossing off'' certain nodes.  Let
$\beta_1=\alpha_1$ and $\beta_2=\alpha_2+\dots+\alpha_n$.  Let $K$ be
the group generated by $T$, $U_{\pm\beta_1}$, and $U_{\pm\beta_2}$. 
Then $K$ is of type $B_2$.  The subsystem $\rtsys(H)$ may contain
either $\beta_1$ or $\beta_2$, but not both.  Similarly $\rtsys(H)$
may not contain $i\beta_1 + j\beta_2$ with $i$ and $j$ positive
integers.  This shows that $K\cap H$ is either $T$ or of type
$A_1$.\medskip

Part (v).  We have that $G=F_4$.  To construct all short $A_2$
complements, take $\phi$ equal to all the short roots in
$\rtsys(F_4)$, so $\phi(H)$ equals all the short roots in $\rtsys(H)$. 
Observe that $\phi$ is of type $D_4$.  By examining each possibility
for $H$ it is easy to verify that $\phi(H)$ is submaximal in a $D_4$
root system and we are done by Corollary
\ref{corollary_sub-maximal_implies_non_spherical}.
%
To construct the long $A_2$ complements, one proceeds similarly with
$\phi$ equal to all the long roots in $\rtsys(F_4)$.
\end{proofof}

\end{document}